%% file: paper.tex
\documentclass[11pt,a4paper]{article}
\usepackage[utf8]{inputenc}
\usepackage{microtype}
\usepackage{amsmath,amsthm,amssymb}
\usepackage{mathrsfs} 
\usepackage[T1]{fontenc}
\usepackage{ae,aecompl}
\usepackage{times}
\usepackage[british]{babel}
\usepackage{aliascnt}
\usepackage[colorlinks,pagebackref,hypertexnames=false]{hyperref} 
\usepackage{url}
\usepackage[alphabetic,backrefs]{amsrefs}
\usepackage{nth}
\usepackage{float}
\usepackage[all]{xy}
\usepackage{tikz-cd} 
\usepackage{bookmark}
\usepackage[margin=1.9cm]{geometry}
\usepackage{enumerate}
\input{init}

\newcommand{\G}{{\mathrm{G}}}
\newcommand{\Id}{{\mathrm{Id}}}

\usepackage{todonotes}
\usepackage{comment}

\allowdisplaybreaks

\usepackage{titling}
\usepackage{changepage}

\begin{document}
	\title{On the quaternionic Goldberg conjecture}
	\author{
	Udhav Fowdar \\
  {University of Warsaw}\\
  u.fowdar@uw.edu.pl}
  \date{}

\maketitle
	
\vspace{-0.5cm}	
\begin{abstract}
We prove a quaternionic K\"ahler analogue of the Goldberg conjecture, namely:\\

Let $(M^8,g,\Om)$ denote a compact $8$-manifold endowed with an $\Sp(2)\Sp(1)$-structure defined by a closed $4$-form $\Omega$. If the metric $g$ is Einstein with non-negative scalar curvature then it is quaternionic K\"ahler i.e. $\Omega$ is parallel. In fact we prove a slightly stronger result allowing for the scalar curvature to be negative under a suitable pinching condition by the square norm of the intrinsic torsion. 
\end{abstract}
	
\begin{adjustwidth}{0.95cm}{0.95cm}
    \tableofcontents
\end{adjustwidth}

\section{Introduction}
Let $(M^{2n},g,J,\omega)$ denote a {compact} almost Hermitian manifold and $\nabla$ its Levi-Civita connection. The Goldberg conjecture \cite{Goldberg1969} can be stated as follows:\medskip

{\indent{Given a $\mathrm{U}(n)$-structure on $M$ determined by $(g,J,\omega)$ such that $d\omega=0$ and $\mathrm{Ric}(g)=\lambda g$, then $\nabla \omega=0$.}}\medskip

\noindent In other words, the conjecture says that the closedness of K\"ahler 2-form together with the Einstein condition imply that the \textit{almost} complex structure $J$ is in fact integrable. As formulated above, the Goldberg conjecture can be extended in a natural way to more general $H$-structures, where $H\subset \mathrm{O}(n)$. More precisely, suppose that there is an underlying differential form $\xi$ associated to an $H$-structure such that $\xi$ is closed and the metric is Einstein. The generalised Goldberg conjecture then asks if this implies that $\nabla \xi=0$. The original Goldberg conjecture corresponds to the case $H=\mathrm{U}({n})$ and $\xi=\omega$. 

The generalised Goldberg conjecture can, in particular, be investigated for the groups $H$ belonging to the Berger's list of Riemannian holonomy groups \cite{Berger}: 
\[
\mathrm{U}(n), 
\qquad\mathrm{SU}(n), \qquad\mathrm{Sp}(n), \qquad\mathrm{Sp}(n)\mathrm{Sp}(1), 
\qquad \mathrm{G}_2, 
\qquad \mathrm{Spin}(7).
\]
In the $\mathrm{U}(n)$ case, the Goldberg conjecture is known to be true in all dimensions if the scalar curvature is non-negative, i.e. $\lambda\geq 0$, \cites{Sekigawa1985, Sekigawa1987}. While the general Goldberg conjecture is still open in the case of negative scalar curvature, there are partial results under additional assumptions, see for instance \cite{OguroSekigawa2001}.
Note also that the hypothesis that $M$ is compact is actually necessary here since there exist non-compact  counterexamples \cites{ApostolovDraghiciMoroianu2001, NurowskiPrzanowski1999}. The odd dimensional Sasakian analogue of the Goldberg conjecture is also known to be true in some cases, see \cites{ApostolovDraghiciMoroianu2006, BoyerGalicki2001}.

Let us now consider the case when $H=\mathrm{G}_2$ and $\xi=\varphi$ corresponds to the $\G_2$ $3$-form. Under the assumption that $\varphi$ is closed, it is known that the scalar curvature is non-positive \cite{Bryant06someremarks}. If the underlying 7-manifold $M$ is compact then Cleyton-Ivanov showed that the Einstein condition indeed implies $\nabla \varphi =0$ \cite{CleytonIvanov2007}; a much shorter proof was subsequently obtained in \cite{Bryant06someremarks}. More recently, it was shown that the conjecture is also true in the non-compact case when the underlying metric is complete \cite{Payne2025}. Thus, in the $\G_2$ case the generalised Goldberg conjecture holds under greater generality than the $\mathrm{U}(n)$ case. Note that one can also formulate the problem for the $\G_2$ $4$-form $\xi=*\varphi$ instead; however, the analogous conjecture is false since nearly parallel $\G_2$-structures provide counterexamples which are Einstein but not torsion free.

In the case of $H=\mathrm{Spin}(7)$, if we ask that the underlying $\mathrm{Spin}(7)$ $4$-form $\xi=\Phi$ is closed then the result of Fern\'andez  \cite{Fernandez1986} implies that this is equivalent to $\nabla \Phi=0$. So the conjecture is automatically true in this case without any assumption of compactness or completeness; the result is purely local, and moreover the Einstein condition is already a consequence of the closedness of $\Phi$. Similarly, for the $H=\mathrm{SU}(n)$ and $\mathrm{Sp}(n)$ cases, if we ask that all the defining differential forms are closed then we automatically obtain that they are parallel.

Lastly we consider the case $H=\mathrm{Sp}(n)\mathrm{Sp}(1)$ for $n\geq 2$, and  with $\xi=\Omega$ the quaternion-K\"ahler 4-form. In \cite{Swann1991}, Swann showed that if $n>2$ then $d\Omega=0$ is equivalent to $\nabla\Omega=0$. Thus, as in the $\Spin(7)$ case, we automatically have that the conjecture is true for $n\geq 3$. This leaves open the remaining case $n=2$. In \cite{Salamon2001} Salamon constructed the first example of a compact $\mathrm{Sp}(2)\mathrm{Sp}(1)$-structure for which $d\Omega=0$ but $\nabla \Omega\neq0$. Subsequently, his construction was generalised and many more examples were obtained in \cite{ContiMadsen2015}. All these examples have closed $\Omega$ but are not Einstein. Thus, this leads to the following quaternionic Goldberg conjecture:\medskip

Given an $\mathrm{Sp}(2)\mathrm{Sp}(1)$-structure on an $8$-manifold $M$ with closed quaternion-K\"ahler $4$-form $\Omega$ and $\mathrm{Ric}(g)=\lambda g$, then $\nabla \Omega=0$.\medskip

\noindent Our main result is a quaternionic analogue of Sekigawa's result in \cites{Sekigawa1985, Sekigawa1987}:
\begin{theorem}\label{main thm}
The quaternionic Goldberg conjecture is true for $M$ compact and $\lambda\geq 0$.
\end{theorem}
In fact our result is slightly stronger. The quaternionic Goldberg conjecture holds under the weaker hypothesis that at each point $\lambda > -\frac{7}{24}|T|^2$, whenever the intrinsic torsion of the $\Sp(2)\Sp(1)$-structure $T$ is non-zero. In other words, the scalar curvature can be negative provided it is suitably pinched. In view of this, it seems reasonable to try to adapt the approach here to extend the result in \cite{Sekigawa1987} to allow for a suitable pinching condition on the scalar curvature in the almost K\"ahler case as well. \\

\noindent\textbf{Strategy of proof.} The overall idea of the proof is quite straightforward: we need to find an integral identity which shows, under the hypothesis of closed $\Omega$ and non-negative Einstein, that the intrinsic torsion vanishes. This is the method employed in the compact $\G_2$ case \cites{CleytonIvanov2007, Bryant06someremarks}, the almost K\"ahler case \cites{Sekigawa1985, Sekigawa1987}, as well as the odd dimensional $K$-contact case in \cites{ApostolovDraghiciMoroianu2006, BoyerGalicki2001}. It is worth noting that there does not appear to be a standard integral identity that can be applied universally to all these $H$-structures, so the derivation does depend on the particular $H$-structure.  Our approach here is mostly based on the methods used in \cite{CleytonIvanov2007} and \cite{Bryant06someremarks} in the $\G_2$ case. Following \cite{CleytonIvanov2007}*{\S 6}, also \cite{Sekigawa1987}, the integral identities are obtained by considering the pairing of the Pontryagin class with the closed quaternion-K\"ahler 4-form, and using the fact that the Levi-Civita connection of an Einstein metric corresponds to a critical point of the Yang-Mills functional. The second step is to rewrite the integrand in a suitable way as a sum of non-negative terms. To do this, we adapt the representation theory method employed by Bryant in \cite{Bryant06someremarks}, see also \cite{fino2026someremarks}, to classify the relevant differential invariants of the underlying $\Sp(2)\Sp(1)$-structures. 

The organisation of the paper is as follows. In section 2, we collect basic facts about $\mathrm{Sp}(2)\mathrm{Sp}(1)$-structures. In section 3, we obtain two integral identities. To obtain the correct integrand we take a linear combination of these. In section 4, we derive the main pointwise identity needed to rewrite the integrand as a sum of non-negative quantities and a divergence part.
Integration with the assumption of $\lambda\geq0$ then yields the desired result.
\\

\noindent\textbf{Acknowledgement.}
I would like to thank Simon Salamon for bringing this problem to my attention several years ago. I am also grateful to Anna Fino for several comments on an earlier draft that helped improved this article.\\

\noindent\textbf{AI disclosure.} 
I had access to OpenAI models through the ChatGPT for Academic Researchers program. The main use of AI in this work, aside from proofreading, was to aid in the search in \S \ref{sec: remaining}. It was not at all obvious to me at the start that the representation theoretic approach would yield an identity of the form \eqref{key identity}, so I used Astra to explore this possibility. This helped establish the feasibility of the approach and substantially reduced the time spent on exploratory calculations and guesswork. Subsequently I independently verified the formula. The manuscript was written entirely by me, and I take full responsibility for its content.

\section{Background on \texorpdfstring{$\mathrm{Sp}(2)\mathrm{Sp}(1)$}{}-structures}

We recall the basics of $\mathrm{Sp}(2)\mathrm{Sp}(1)$-structures, and fix our notation and conventions for the rest of this article. We refer the reader to \cites{Salamon1989, Salamon1982} for further details. 

An $\mathrm{Sp}(2)\mathrm{Sp}(1)$-structure on an $8$-manifold $M$ is determined by a $4$-form $\Omega$ which can locally be expressed as:
\begin{equation}
    \Omega=\frac{1}{2}(\omega_1 \wedge\omega_1+\omega_2 \wedge\omega_2+\omega_3 \wedge \omega_3),
\end{equation}
where
\begin{align*}
\omega_1 &= e^{12} + e^{34} + e^{56} +e^{78},\\ 
\omega_2 &= e^{13} + e^{42} + e^{57} +e^{86},\\
\omega_3 &= e^{14} + e^{23} + e^{58} +e^{67},
\end{align*}
and $\{e^i\}_{i=1}^8$ denotes a local co-framing on $M$. Such a co-framing is said to be $\Sp(2)\Sp(1)$ adapted. At each point the stabiliser of $\Omega$ in $\mathrm{GL}(8,\R)$ is isomorphic to $\Sp(2)\Sp(1) \subset \SO(8)$ and hence $\Omega$ defines a Riemannian metric $g$ and volume form $\vol= \frac{1}{30}\Omega\wedge\Omega$. We denote by $\nabla^{LC}$ the Levi-Civita connection of $g$ and by $*$ the Hodge star operator. 

\begin{remark}
It is worth noting that although the above $2$-forms $\om_1,\om_2,\om_3$ need not exist globally on $M$, they span a globally well-defined rank $3$ subbundle of $\Lm^2$; we denote this by $\Lambda^2_3$. Since one can choose local trivialisation of this bundle, we can always locally reduce the structure group from $\Sp(2)\Sp(1)$ to $\Sp(2)$.   
\end{remark}

For simplicity we shall write $\mathfrak{h}$ for the Lie algebra $\mathfrak{sp}(2)\oplus\mathfrak{sp}(1)$ and denote by $\mathfrak{h}^\perp$ its orthogonal complement in $\mathfrak{so(8)}\cong \Lambda^2$. Throughout the article, we shall use the metric to identify 1-forms with vector fields, and hence also 2-tensors with endomorphisms. We define the canonical connection $\nabla^c$ by projection of $\nabla^{LC}$ into $\mathfrak{h}$. More precisely, we have
\begin{equation}
\nabla^{LC}_XY-\nabla^c_XY=T(X)(Y), \label{def of canonical connection}    
\end{equation}
where $T \in \Lambda^1\otimes \mathfrak{h}^{\perp}$ is the intrinsic torsion of the $\Sp(2)\Sp(1)$-structure. We shall use the lowering index convention $T_{i,k}^{\  j}=T_{i,jk}$ to view $T$ as a $2$-form valued $1$-form. By the definition of $\nabla^c$, it follows that $\nabla^c\Omega=0$. The space of intrinsic torsion has dimension $\dim(\Lm^1 \otimes \mathfrak{h}^\perp) = 8 \times (28 - 13) =120$. To describe the irreducible components of $T$, we first need to recall some representation theory.

Following Salamon's $E$-$H$ formalism in \cite{Salamon1982}, we denote by $E\cong \C^4$ and $H\cong \C^2$ the standard complex representations of $\Sp(2)$ and $\Sp(1)$ respectively. In this notation the complexified space of $1$-forms can be identified with $\Lambda^1_{\C}\cong E\otimes H$. We also have $\mathfrak{h_{\C}}=S^2 E\oplus S^2H$ and $\mathfrak{h}^\perp_{\C}= \Lambda_0^2E\otimes S^2H$. To simplify the notation we shall henceforth omit the complexification sign and it should be understood that we are considering the real subspaces underlying the complex representations of $\mathrm{Sp}(2)\Sp(1)$ whenever a module can be realised as a subspace of $(\otimes^pE)\otimes (\otimes^qH)$ for $p+q$ even.

The space of $2$-forms splits into irreducible modules as:
\[
\Lambda^2=\Lambda^2_3 \oplus \Lambda^2_{10}\oplus \Lambda^2_{15},
\]
where
\begin{equation}\label{decomposition of 2 forms}
\begin{aligned}
    S^2H \cong \Lambda^2_3
    &= \{\alpha\in\Lambda^2 \mid *(\alpha\wedge\Omega)=+5\alpha\}, \\
    S^2E \cong \Lambda^2_{10}
    &= \{\alpha\in\Lambda^2 \mid *(\alpha\wedge\Omega)=-3\alpha\}, \\
    \Lambda^2_0E\otimes S^2H \cong \Lambda^2_{15}
    &= \{\alpha\in\Lambda^2 \mid *(\alpha\wedge\Omega)=+\alpha\}.
\end{aligned}
\end{equation}
Here the notation $\Lambda^p_k$ should be understood as saying that this is a subspace of $\Lambda^p$ of real dimension $k$. We denote by $P_k$ the projection maps
$$P_{k}:\Lambda^2 \to \Lambda^2_k,$$ where $k=3,10,15$. 
Similarly, for the $\Lambda^3$ we have
\begin{equation}
\Lambda^3=  \Lambda^3_8 \oplus \Lambda^3_{16}\oplus \Lambda^3_{32}\cong (E \otimes H)\oplus (E \otimes S^3 H) \oplus (V_{1,1}\otimes H ),\label{decom of lm3}
\end{equation}
where $V_{1,1}$ denotes the irreducible complex representation of $\Sp(2)$ defined by $E \otimes \Lm^2_0E \cong V_{1,1}\oplus E$. The space $V_{1,1}$ has complex dimension 16 and $\Lm^3_{32}$ corresponds to the real subspace underlying $V_{1,1}\otimes H$. Our convention here is $E\cong V_{0,1}$ and $\Lambda^2_0E \cong V_{1,0}$. 
Lastly we recall the decomposition of the space of symmetric 2-tensors:
\begin{equation}\label{decom symmetric tensor}
    S^2(\Lambda^1)\cong \R \oplus \Lambda^2_0E \oplus S^2H\otimes S^2E. 
\end{equation}
In particular, this shows that the traceless Ricci tensor has 2 irreducible components of dimensions 5 and 30 respectively; the 1-dimensional component is of course just the scalar curvature. 

We can now identify the four irreducible components of the intrinsic torsion as follows: 
\begin{equation}
\nabla^{LC} \Omega \simeq T \in \Lambda^1 \otimes \mathfrak{h}^\perp \cong (E\otimes H) \oplus (E\otimes S^3H)\oplus (V_{1,1}\otimes H) \oplus (V_{1,1}\otimes S^3H).\label{torsion decom}
\end{equation}
The result of Swann in \cite{Swann1991} says that the first 3 summands in \eqref{torsion decom} can be identified with the 3 components of $d\Omega \in \Lambda^5 \cong \Lambda^3$ in view of \eqref{decom of lm3}. Thus, if $d\Omega=0$ then $\nabla^{LC} \Omega\simeq T \in V_{1,1}\otimes S^3H$. 
More explicitly, one can show:
\begin{proposition}\label{prop; torsion}
The condition $d\Omega=0$ is equivalent to 
\begin{equation}
\begin{matrix}
d\om_1 &=& 0 &+& \alpha_{12} \w \om_2 &+& \alpha_{13} \w \om_3  &+& \beta_1,\\
d\om_2 &=& - \alpha_{12} \w \om_1 &+& 0  &+& \alpha_{23} \w \om_3  &+& \beta_2,\\
d\om_3 &=& - \alpha_{13} \w \om_1 &-& \alpha_{23} \w \om_2 &+& 0 &+&\beta_3,
\end{matrix}
\end{equation}
where $\alpha_{ij}\in \Lambda^1$ and $\beta_i \in \Lambda^3_{32}$, and furthermore $\beta_1 \w \om_1 + \beta_2 \w \om_2 + \beta_3 \w \om_3 =0 $. In other words, the intrinsic torsion of a closed quaternion K\"ahler $4$-form $\Om$ can be identified with the tensor:
    \[ \beta_1 \otimes \om_1 + \beta_2 \otimes \om_2 + \beta_3 \otimes \om_3 \in V_{1,1}\otimes (H \oplus S^3 H)\]
    subject to the relation
    $ \beta_1 \w \om_1 + \beta_2 \w \om_2 + \beta_3 \w \om_3=0 \in V_{1,1}\otimes H$.
\end{proposition}
To reconcile these two descriptions of the intrinsic torsion, we first recall that from \eqref{def of canonical connection} we have
\begin{equation}
\nabla^{LC}_X \Omega =  T(X)\diamond \Omega,    
\end{equation}
where for a $2$-tensor $\alpha$ the endomorphism action $\diamond$ is defined by
\[
\alpha\diamond \Omega=
\sum_{i,j=1}^8   \alpha_{ij}e^i \wedge (e_j\ip \Omega).
\]
When $d\Omega=0$, the two descriptions of the intrinsic torsion given above are explicitly related by:
\begin{equation}\label{T in terms of beta}
T=\frac{1}{2} \sum_{i=1}^8\sum_{j=1}^3
(e_i \ip \om_j)\otimes (e_i \ip \beta_j).
\end{equation}
\begin{remark}
Note that $\dim(V_{1,1}\otimes S^3H)=64$ while $\dim(\Lambda^1\otimes \mathfrak{h}^\perp)=120$, so the closedness condition of $\Omega$ imposes only $56$ equations. As a comparison with $\G_2$-structures, the closedness condition of $\varphi$ imposes $35$ equations leaving only a $14$-dimensional space of possible intrinsic torsions. Heuristically, this suggests that the hypothesis of generalised Goldberg conjecture are somewhat stronger in the $\G_2$ than in the $\Sp(2)\Sp(1)$ case.
\end{remark}
Given an affine connection $\nabla$, we define its curvature by $R(X,Y)Z=\nabla_X(\nabla_YZ)-\nabla_Y(\nabla_XZ)-\nabla_{[X,Y]}Z$ and we follow the convention $R_{ijkl}= g(R(e_i,e_j)e_l,e_k)$. In particular, for $\nabla^{LC}$ we have $\Ric_{ij}=R_{ikjk}^{LC}$. 
We define the model quaternion-K\"ahler curvature tensor by 
\begin{equation}
R_0=4\sum_{i=1}^{3}u_i \otimes u_i+2\sum_{i=1}^{10}v_i\otimes v_i \in S^2(\Lambda^2_{3})\oplus S^2(\Lambda^2_{10}),
\end{equation}
where $u_i$ and $v_i$ denote an orthonormal basis for $\Lambda^2_3$ and $\Lambda^2_{10}$ respectively. It is easy to see that $R_0$ corresponds to the endomorphism $4P_{3}+2P_{10}$ using $g$, and up to an overall factor, it coincides with the Riemann curvature tensor of $\mathbb{HP}^2$. In this normalisation,  it is easy to check that  $(R_0)_{ikjk}=4g_{ij}$.

\begin{remark}[Inner product conventions]
To avoid having to distinguished between the differential form  inner product 
and the tensor inner product, we shall only use the full
index contraction convention i.e. in an orthonormal frame,
\[
g(A,B)=A_{i_1\cdots i_r}B_{i_1\cdots i_r},
\qquad \text{and} \qquad
|A|^2=A_{i_1\cdots i_r}A_{i_1\cdots i_r}.
\]
So for instance, in this convention we have $|\omega_i|^2=8$.
\end{remark}

In \cite{Witt}, it is shown that when $d\Omega=0$, the component of the Ricci tensor in $\Lambda^2_0E$, in view of \eqref{decom symmetric tensor}, vanishes. Using the $\diamond$ operator, the method used in \cite{Bryant06someremarks} in the $\G_2$ case, see also \cite{fino2026someremarks}, can be adapted to our $\Sp(2)\Sp(1)$ situation to show that when $d\Omega=0$, the traceless Ricci curvature is given explicitly by:
\begin{align}
    \Ric_0(g) \diamond \Omega = \pi^{4-}_{30}\Big(
    &-4(*( \pi^{4-}_{30}(d\beta_1 +\alpha_{21}\w \beta_2 + \alpha_{31} \w \beta_3) \w \om_2  ) \w \om_3)\label{ricci formula}\\
    &-4(*( \pi^{4-}_{30}(d\beta_2 +\alpha_{32}\w \beta_3 + \alpha_{12} \w \beta_1) \w \om_3  ) \w \om_1)\nonumber\\
    &-4(*( \pi^{4-}_{30}(d\beta_3 +\alpha_{13}\w \beta_1 + \alpha_{23} \w \beta_2) \w \om_1  ) \w \om_2)\nonumber\\ 
    &+
    4( 
    *(\beta_{1} \w \beta_{2}) \w \om_3 +
    *(\beta_{2} \w \beta_{3}) \w \om_1 +
    *(\beta_{3} \w \beta_{1}) \w \om_2 
    )
    \Big),\nonumber
\end{align}
where $\pi^{4-}_{30}$ denotes the projection map $\Lambda^4\to \Lambda^4_{30} \subset \Lambda^{4-}\cong \Lambda^2_0E \oplus S^2E \otimes S^2 H=\Lambda^4_{5}\oplus \Lambda^4_{30}$, where $\Lambda^{4-}$ denotes the space of anti-self-dual $4$-forms.

From now on we shall always assume that $(M,g,\Omega)$ is a compact Einstein 8-manifold with $d\Omega=0$. The positive scalar curvature assumption will only be needed in last step of the proof.

\section{The integral identities}

\subsection{The first integral identity.} 
The first integral identity will be obtained using the fact that the Levi-Civita connection of an Einstein manifold corresponds to a critical point of the Yang-Mills functional \cite{Besse2008}*{\S 16}. More precisely, we shall consider a variation of the Levi-Civita connection using the intrinsic torsion. The resulting variation formula for the Yang-Mills functional will then yield the desired identity involving only second order invariants of the $\Sp(2)\Sp(1)$-structure. 
We now give the details. 

As before we fix an adapted local orthonormal co-framing $\{e^i\}$ so that the Levi-Civita connection $\nabla^{LC}$ can be identified with the matrix valued $1$-form $A_{LC} \in \Lambda^1\otimes \mathfrak{so}(8)$ and the canonical connection $\nabla^c$ with $A_c \in \Lambda^1\otimes \mathfrak{h}$. Hence we have
\[
A_{LC}= A_c+T.
\]
To be precise, the convention here is $A^i_j=g(\nabla_{e_k}e_j,e_i)e^k$ and $T^i_j=g(T(e_k)e_j,e_i)e^k$. 
Considering the family of metric connections defined by $A_t:=A_c+tT$, we can compute the associated curvature 2-forms as
\begin{align}
    F_t:=&\ dA_t+A_t \wedge A_t\nonumber\\
    =&\ (dA_c+A_c\wedge A_c)+ t (dT+A_c \wedge T+T\wedge A_c) + t^2 (T \wedge T)\nonumber\\
    =&\ F_c + t d^{c}T + t^2(T \wedge T), \label{equ: FAt}
\end{align}
where $F_c$ denotes the curvature of $A_c$ and $d^cT$ denotes the exterior covariant derivative of $T$ with respect to $\nabla^c$.  In indices, $F^k_l=\frac{1}{2}R_{ijkl}e^i \wedge e^j$. It follows from the above that $$F_c \in \Lambda^2 \otimes \mathfrak{h},\qquad  d^cT\in \Lambda^2 \otimes \mathfrak{h}^{\perp}, \qquad B:=T\wedge T\in \Lambda^2 \otimes \mathfrak{h}.$$ 
The latter is just consequence of the fact that $[\mathfrak{h}^{\perp} , \mathfrak{h}^{\perp}]\subset \mathfrak{h}$, see for instance \cite{FowdarSaEarp2024}*{Lemma 2.3}. Viewing the curvature as a section of $\Lambda^2 \otimes \Lambda^2$, we have
\begin{align}
    d^cT &= (\Id \otimes P_{15})(F_{1})\label{equ: 2222}\\
        F_c+B &= (\Id \otimes (P_{3}+ P_{10}))(F_{1}),\label{equ: 3333}
\end{align}
where $\Id \otimes P_{k}$ denotes the obvious projection map with $\Id$ acting the first two indices and $P_k$ on the last two indices, and $F_1$ corresponds to the Riemannian curvature tensor $R^{LC}$. 

Observe also that the exterior covariant derivative of $T$ with respect to $\nabla^{LC}$ is given by 
\begin{equation}\label{equ: dLCT}
d^{LC}T:=dT+A_{LC}\wedge T+T\wedge A_{LC} = d^cT+ 2 T\wedge T,
\end{equation}
which in index notation is expressed as $(d^{LC}T)^k_l=\frac{1}{2}(\nabla_i^{LC}T_{j,kl}-\nabla_{j}^{LC}T_{i,kl})e^{i}\wedge e^j$.

We now consider the Yang-Mills functional for the above 1-parameter family of metric connections:
\[
E(t)= \frac{1}{4}\int_M (F_{t})_{ijkl}(F_{t})_{ijkl}\vol.
\]
Taking the time derivative at $t=1$ we get
\begin{align*}
    E'(1) = \int_M \nabla^{LC}_{i}T_{j,kl}R^{LC}_{ijkl} \vol
    = -\int_M T_{j,kl}\nabla^{LC}_iR_{ijkl}^{LC}\vol
    = \int_M T_{j,kl}(\nabla^{LC}_kR_{ijli}^{LC}+\nabla^{LC}_lR_{ijik}^{LC})\vol
    =0,
\end{align*}
where the first equality follows \eqref{equ: FAt} and \eqref{equ: dLCT}, the second one by integration by parts, the third one by the second Bianchi identity, and the last one uses the Einstein condition i.e. $R_{ijik}=\Ric_{jk}=\lambda g_{jk}$. 

This is of course a well-known fact since the metric $g$ is Einstein and hence its Levi-Civita connection is a critical point of the Yang-Mills functional, see \cite{Besse2008}*{\S 16.34}. 

The upshot now is that using \eqref{equ: FAt} we can rewrite this as
\[
0=\int_Mg(F_1,d^cT+2B) \vol = \int_Mg(F_c+d^cT+ B,d^cT+2B)\vol =
\int_M\big(g(F_c+B,2B)+ |d^cT|^2\big)\vol,
\]
where we used that $\mathfrak{h}$ and $\mathfrak{h}^\perp$ are orthogonal. Thus, we obtain the first integral identity:
\begin{equation}\label{first integral identity}
 I_1:= \int_M\Big(g\big((\Id \otimes (P_{3}+ P_{10}))R^{LC},2B\big)+ |d^cT|^2\Big)\vol  =0.
\end{equation}

\subsection{The second integral identity.}
The second integral identity will be obtained using Chern-Weil theory following an analogous strategy used by Cleyton-Ivanov in \cite{CleytonIvanov2007}*{\S 6}. More precisely, 
given any 2 connections $A_1$ and $A_2$ on the tangent bundle we have
\[
 \mathrm{tr}(F_{A_1}\wedge F_{A_1}) - \mathrm{tr}(F_{A_2}\wedge F_{A_2}) = d\alpha,
\]
for some $3$-form $\alpha$, and the cohomology class $[ \mathrm{tr}(F_{A_1}\wedge F_{A_1})]$ corresponds to the first Pontryagin class of $M$ up to an overall constant. Thus, we have
\[
\int_M \mathrm{tr}(F_{A_1}\wedge F_{A_1}) \wedge \Omega - \int_M \mathrm{tr}(F_{A_2}\wedge F_{A_2}) \wedge \Omega=\int_M d(\alpha\wedge \Omega)=0,
\]
where we used that $d\Omega=0$ and Stokes' theorem. Thus, the integral of 
$\mathrm{tr}(F_{A_t}\wedge F_{A_t}) \wedge \Omega$ is independent of the parameter $t$ for the family of connections $A_t=A_c+tT$. Extracting the term with the factor of $t^2$ gives:
\begin{equation}\label{equ: 111}
    \int_M \mathrm{tr}(F_c \wedge B+  B \wedge F_c+d^cT \wedge d^cT) \wedge \Omega=0.
\end{equation}
Next we want to rewrite \eqref{equ: 111} in a similar way as \eqref{first integral identity} by decomposing into $\mathfrak{h}$ and $\mathfrak{h}^\perp$ components. To this end, we first define the $\Sp(2)\Sp(1)$-equivariant map $L:\Lambda^2 \to \Lambda^2$ by
\begin{equation}\label{definition of L}
L(\alpha) = *(\alpha \wedge \Omega).
\end{equation}
The irreducible summands $\Lambda^2_k$ in \eqref{decomposition of 2 forms} correspond precisely to the eigenspaces of $L$. Using \eqref{definition of L}, we can rewrite \eqref{equ: 111} equivalently as
\begin{equation}\label{eq: 44444}
    \int_M \Big( g( (L\otimes \Id) F_c,2B)+ g((L\otimes \Id) d^cT,d^cT)\Big)\vol = 0,
\end{equation}
where as before $L\otimes \Id$ should be understood as acting on $\Lambda^2\otimes \Lambda^2$ with $L$ acting on the first copy of $\Lambda^2$ and $\Id$ acting on the second one. 

We now want to replace the canonical curvature $F_c$ in \eqref{eq: 44444} by the Riemannian curvature $R^{LC}$ so that it is of the same form as \eqref{first integral identity}. Applying $L\otimes\Id$ to \eqref{equ: 3333}, we get  $(L\otimes\Id)F_c = (L \otimes (P_{3}+ P_{10}))R^{LC}-(L\otimes\Id)B$, which can now substitute in \eqref{eq: 44444} to get
\begin{equation}\label{second integral identity}
I_2:=    \int_M \Big( g( (L \otimes (P_{3}+ P_{10}))R^{LC},2B)+ g((L\otimes \Id) d^cT,d^cT)\Big)\vol = 0.
\end{equation}
Note that here we also used the fact that $\mathrm{tr}(B\wedge B)=\mathrm{tr}(T \wedge T \wedge T \wedge T)=0$ since $T$ is a matrix of 1-forms.

\subsection{The combination}
Recall that if $T=0$ then 
\[
R^{LC}=R^c \in \mathrm{ker}\big(\wedge:S^2(\mathfrak{h}) \to \Lambda^4\big)\cong \mathbb{R}\oplus S^4E,
\]
where the $\R$ component is generated by $R_0$ with its coefficient determined by the scalar curvature, and $S^4E$ corresponds to space of curvature tensors of hyperK\"ahler manifolds, see for instance \cite{Salamon1989}.

The main issue at this point is that both \eqref{first integral identity} and \eqref{second integral identity} contain part of the curvature component of $R_{LC}$ which lies in $S^4E$. By the above we see that this component of curvature cannot be controlled by $T$; hence we need to eliminate this, since such a term cannot give rise to an integral identity which implies the vanishing of $T$. As $S^4E\subset \Lambda^2_{10}\otimes \Lambda^2_{10}$, it suffices to eliminate the $(P_{10} \otimes P_{10})R^{LC}$ term appearing in the integral identities. 

From \eqref{decomposition of 2 forms} we see that $L$ has eigenvalue $-3$ on $\Lambda^2_{10}$, hence we have  $(L+3\Id)\Lambda^2 \subset \Lambda^2_3\oplus\Lambda^2_{15}$. Comparing \eqref{first integral identity}
and \eqref{second integral identity}, we can thus take the linear combination $3I_1+I_2$ to get:
\begin{equation}\label{main integral identity}
 \int_M\Big(g\big(((L+3\Id) \otimes (P_{3}+ P_{10}))R^{LC},2B\big)+ g\big(((L+3\Id)\otimes \Id) d^cT,d^cT\big)\Big)\vol =0,
\end{equation}
where we also note from \eqref{decomposition of 2 forms} that $L+3\Id=8P_3+4P_{15}$. In particular, the second term in \eqref{main integral identity} is non-negative:
\begin{align*}
g\big(((L+3\Id)\otimes \Id) d^cT,d^cT\big) &= g\big((8P_3\otimes \Id) d^cT,d^cT\big)+g\big((4P_{15}\otimes \Id) d^cT,d^cT\big)\\
&=8|(P_3\otimes \Id) d^cT|^2+4|(P_{15}\otimes \Id) d^cT|^2.
\end{align*}
Since the integrand of \eqref{main integral identity} does not contain the $S^4E$ component of the Riemann curvature tensor, it follows that it is completely determined by $T$, its derivative, and the scalar curvature. 

Recalling the definition of the model quaternion-K\"ahler curvature tensor $R_0$, we can write $R^{LC}= \frac{\lambda}{4} R_0+\hat{R}$. Contracting the latter with the metric, we get $\Ric_{ij}=\lambda g_{ij}+ \hat{R}_{ikjk}$ and hence from the Einstein hypothesis, we deduce $\hat{R}_{ikjk}=0$. The upshot is that we can now rewrite the first term in \eqref{main integral identity} as:
\begin{align*}
    g\big(((L+3\Id) \otimes (P_{3}+ P_{10}))R^{LC},2B\big) &= 2\lambda g\big((P_3 \otimes P_{3})R_0,2B\big) + g\big(((L+3\Id) \otimes (P_{3}+ P_{10}))\hat{R},2B\big) \\
    &= \lambda g\big(\sum_{i=1}^3\omega_i\otimes \omega_i,2B\big) + g\big(((L+3\Id) \otimes (P_{3}+ P_{10}))\hat{R},2B\big) \\
    &= 4\lambda |T|^2 + g\big(((L+3\Id) \otimes (P_{3}+ P_{10}))\hat{R},2B\big).
\end{align*}
where in the first equality we used the definition of $R_0$ and the fact that $L+3\Id=8P_3+4P_{15}$, and for the last equality we use the fact that $g\big(\sum_{i=1}^3\omega_i\otimes \omega_i,2B\big)$ has to be equal to a constant multiple of $|T|^2$ since there is exactly one copy of the trivial representation in $S^2(V_{1,1}\otimes S^3H)$; so by Schur's lemma it suffices to do an algebraic check to obtain the coefficient factor. 

Hence we can now rewrite \eqref{main integral identity} as
\begin{equation}\label{main integral identity 2}
 \int_M\Big(\lambda |T|^2 + \mathcal{H}_{rem}  \Big)\vol =0,
\end{equation}
where 
\begin{equation*}
\mathcal{H}_{rem}:=    \frac{1}{4}g\big(((L+3\Id) \otimes (P_{3}+ P_{10}))\hat{R},2B\big)+ 2|(P_3\otimes \Id) d^cT|^2+|(P_{15}\otimes \Id) d^cT|^2 .
\end{equation*}
The last step of the proof is to show that the integral of remainder term $\mathcal{H}_{rem}$ is non-negative. More precisely, we shall show that $\mathcal{H}_{rem}$ can in fact be expressed (pointwise) as a sum of a divergence term and a non-negative term.

\section{The pointwise algebraic relation}

Observe that $\mathcal{H}_{rem}$ involves at most quadratic terms which are second order in the $\Sp(2)\Sp(1)$-structure. More precisely, it involves terms of the form: $\nabla^cT*\nabla^cT$,  $R^{LC}*T*T$, $\nabla^cT*T*T$, and $T*T*T*T$. As we have already removed the $S^4E$ component of $R^{LC}$ in \eqref{main integral identity}, the component of $R^{LC}$ appearing in $\mathcal{H}_{rem}$ is in fact expressible in terms of $\nabla^cT$ and $T*T$ as well. So the $\mathcal{H}_{rem}$ only involves terms of the form: $\nabla^cT*\nabla^cT$, $\nabla^cT*T*T$ and $T*T*T*T$, where by hypothesis we have $T \in V_{1,1}\otimes S^3H$. Thus, we are only interested in the trivial representations occurring in the modules:
\begin{equation}
S^2(E\otimes H \otimes  V_{1,1} \otimes S^3H),\qquad 
E\otimes H \otimes  V_{1,1} \otimes S^3H\otimes S^2( V_{1,1}\otimes S^3H), \qquad
S^4(V_{1,1}\otimes S^3H).
\end{equation}
Furthermore, we also need to consider the divergence terms since these do not contribute to the integral. The above discussion suggests that we should consider divergence terms of the form
$\nabla_i^{LC} \alpha_i$, where $\alpha$ is a 1-form {cubic} in the torsion i.e. $\alpha \simeq T*T*T$. So we also need to know the $E\otimes H$ modules occurring in:
\begin{equation*}
      S^3(V_{1,1}\otimes S^3H).
\end{equation*}
Hence we have reduced the problem of expressing $\mathcal{H}_{rem}$ as a sum of a divergence term and a non-negative term to a finite (albeit long) representation theoretic calculation. 

\subsection{Divergence term}
We consider the possible $1$-forms which are cubic in torsion and hence give rise to suitable divergence term. 
A straightforward representation theoretic calculation gives:
\begin{align*}
  S^3(V_{1,1}\otimes S^3H) &\cong S^3(V_{1,1})\otimes S^3(S^3H)\oplus S_{(2,1)}V_{1,1}\otimes S_{(2,1)}(S^3H) \oplus \Lambda^3 V_{1,1}\otimes \Lambda^3 (S^3H)  \\
  &\cong 2 E \otimes H \oplus \bigoplus_{i\in J_0} R_i,
 \end{align*}
where none of the irreducible summands $R_i$ are isomorphic to $E \otimes H$. Thus, it follows that there is a 2 dimensional space of possible $1$-forms. We can construct a basis as follows. First we define two symmetric tensors which are quadratic in torsion by: 
\begin{equation}\label{ioioioi}
{Y}_{ab}:=T_{a,jk}T_{b,jk}, \qquad Z_{ab}:=T_{i,aj}T_{i,bj}. 
\end{equation}
Contracting with $T$ again gives the desired basis:
\[
(\alpha_1)_i:=Y_{ab}T_{a,ib}, \qquad (\alpha_2)_i:={Z}_{ab}T_{a,ib} 
\]
It is an easy algebraic calculation to check that $\alpha_1$ and $\alpha_2$ are indeed independent. So we can consider divergence term of the form $a\nabla_i(\alpha_1)_i+b\nabla_i(\alpha_2)_i$, where $a,b\in \R$.
\subsection{The remaining terms}\label{sec: remaining}

The divergence terms we found above correspond to two generators of the trivial representations in:
\begin{gather}\label{qoqoqoq}
\nabla^cT * T*T\in E\otimes H \otimes  V_{1,1} \otimes S^3H\otimes S^2( V_{1,1}\otimes S^3H)\cong 6 \R \oplus \bigoplus_{i \in J_1} R_i 
\end{gather}
where the irreducible representations $R_i$ are all non-trivial. A similar representation theoretic calculation gives:
\begin{gather}
\nabla^cT * \nabla^cT\in S^2(E\otimes H \otimes  V_{1,1} \otimes S^3H) \cong 8 \R \oplus \bigoplus_{i \in J_2} R_i, \label{equ: ppopo}\\
T*T* T*T\in S^4(V_{1,1}\otimes S^3H) \cong 7 \R \oplus \bigoplus_{i \in J_3} R_i.\label{quartic term}
\end{gather}
Thus, there is up to a total of $21=(6+8+7)$ terms that can be used to rewrite $\mathcal{H}_{rem}$. We have already found 2 of these terms, namely $\mathrm{div}(\alpha_1)$ and $\mathrm{div}(\alpha_2)$. This leaves up to $19$ terms to be determined.

However, we also need to account for the algebraic relations coming from the first Bianchi identity. More precisely, recall that:
\begin{equation}
    d^cT = (\Id \otimes P_{15})R^{LC} \quad \text{and} \quad R^{LC} \in \mathcal{R}:= \mathrm{ker}\big(\wedge:S^2(\Lambda^2) \to \Lambda^4\big),
\end{equation}
where
\begin{equation*}
    \mathcal{R}\cong 2\R \oplus S^4E \oplus \Lambda^2_0E \oplus S^2_0(\Lambda^2_0E)\oplus (\Lambda^2_0 E \oplus S^2_0(\Lambda^2_0E))\otimes S^4H \oplus (2S^2E\oplus \Lambda^2_0E \oplus V_{1,2})\otimes S^2H \oplus S^4H.
\end{equation*}
On the other hand, we also have
\begin{equation}\label{equ: nabla T}
\nabla^cT\in E\otimes H\otimes V_{1,1}\otimes S^3H\cong (\Lambda^2_0 E \oplus S^2E \oplus S^2_0(\Lambda^2_0E)\oplus V_{1,2}) \otimes (S^4H\oplus S^2H).
\end{equation}
Thus, the Bianchi identity shows that $d^cT$ has no terms in $S^2_0(\Lambda^2_0E)\otimes S^2H$ and $S^2E \otimes S^4H$ as $R^{LC}$ does not. 
Since the highest order terms of $d^cT$ are determined by $\nabla^cT$, explicitly 
\begin{equation*}
(d^cT)_{ijkl}
=
\nabla^c_iT_{j,kl}-\nabla^c_jT_{i,kl}
+T_{i,jp}T_{p,kl}-T_{j,ip}T_{p,kl},
\end{equation*}
it follows that the corresponding components of $\nabla^cT$ have to vanish at highest order as well; it suffices to check that skew-symmetrisation is non zero on each of these irreducible modules. 
The two representations $S^2_0(\Lambda^2_0E)\otimes S^2H$ and $S^2E \otimes S^4H$ contribute to 2 copies of $\R$ in \eqref{equ: ppopo}, namely it is their squared norms; hence since from the first Bianchi identity these vanish, this reduces the number of possible unknown invariants to $17=19-2$.

Lastly, we impose the Einstein condition.
In \cite{Witt}*{Proposition 33 (ii)}, Witt showed that for closed quaternion-K\"ahler $4$-form, the traceless component of the Ricci curvature lies in $S^2E \otimes S^2H$ i.e. the 5 dimensional part in $\Lambda^2_0E$ vanishes. Thus, we get 
\begin{gather}
    \Ric_0 * \nabla^c T \in S^2E \otimes S^2H\otimes E \otimes H \otimes V_{1,1}\otimes S^3H \cong  \R \oplus \bigoplus_{i} R_i\\
    \Ric_0 * T* T \in S^2E \otimes S^2H \otimes S^2(V_{1,1}\otimes S^3H) \cong 2\R \oplus \bigoplus_{i} R_i 
\end{gather}
Hence among the remaining $17$ terms, $3$ terms vanish from the Einstein condition thereby leaving $14$ only unknown invariants. More precisely, the vanishing of $\Ric_0$ implies that the differential invariant in $S^2E\otimes S^2H$ generated by the (highest order) component of $\nabla^cT$ in \eqref{equ: nabla T} is zero; this can also be deduced explicitly using \eqref{ricci formula}. 

So to summarise, the above argument shows that $\mathcal{H}_{rem}$is expressible as a linear combination of $\mathrm{div}(\alpha_1)$, $\mathrm{div}(\alpha_2)$ and these 14 terms. It thus remains to find these 14 terms and derive the relation. Before doing so, we first explain the strategy of proof. Since we are looking to derive a pointwise identity, we can assume without loss of generality that $\Omega$ and $g$ are given by the standard Euclidean models. We then choose possible tensors for $T \in V_{1,1}\otimes S^3H$ and $\nabla^cT \in E\otimes H \otimes V_{1,1}\otimes S^3H$ subject to the algebraic constraints imposed by the Bianchi identity and vanishing traceless Ricci tensor.  
In fact to find suitable choices for $T$, we first select arbitrary $\beta_1, \beta_2 \in \Lambda^3_{32}$ and set $\beta_3=\omega_1\diamond \beta_2-\omega_2\diamond \beta_1=J_1 \beta_2-J_2 \beta_1$. This guarantees that
$\sum \beta_i\otimes \omega_i \in V_{1,1}\otimes S^3H$, see Proposition \ref{prop; torsion}. The intrinsic torsion $T$ is then obtained using \eqref{T in terms of beta}. Here is a simple concrete  example to illustrate. Set
\begin{equation*}
\beta_1=0,\qquad
\beta_2=-e^{135}-e^{146}+e^{236}-e^{245},\qquad
\beta_3=e^{136}-e^{145}+e^{235}+e^{246}.
\end{equation*}
This gives
\[
T= -e^1\otimes (e^{15}-e^{26})+e^2\otimes (e^{16}+e^{25}),
\]
and we can also set $\nabla^cT=0$ and
\[
\hat{R}= 2(e^{12}+e^{56})\otimes (e^{12}+e^{56})- (e^{15}-e^{26})\otimes (e^{15}-e^{26})-(e^{16}+e^{25})\otimes (e^{16}+e^{25})
\]
It is a straightforward calculation to check that both the Bianchi and Einstein conditions hold for the above choices. All these calculations are purely algebraic and hence are easiest done with the help of a computer algebra software; we used Maple. It suffices to make enough choices for the tensors $T$, $\nabla^cT$, $\hat{R}$ to extract 16 linearly independent equations and hence determine all the coefficients. The final result is \eqref{key identity}. With this strategy in mind, it only remains to find the generators for the 14 remaining terms.

From the previous discussion, there are 5 irreducible representations in \eqref{equ: nabla T} that contain non-zero components of $\nabla^cT$; 2 were removed by the Bianchi identity and 1 was removed by the Einstein condition. The 5 remaining trivial modules in \eqref{equ: ppopo} are generated precisely by their square norm: 
\[
|\Pi_{25}(\nabla^cT)|^2,\quad
|\Pi_{15}(\nabla^cT)|^2,\quad 
|\Pi_{70}(\nabla^cT)|^2,\quad
|\Pi_{175}(\nabla^cT)|^2, \quad 
|\Pi_{105}(\nabla^cT)|^2,
\]
where as usual the inner product is $|\nabla^cT|^2=\nabla^c_{i}T_{j,kl}\nabla^c_{i}T_{j,kl}$ and each $\Pi_k$ denotes the projection map into an irreducible module of dimension $k$; for instance, $\dim(\Lambda^2_0E \otimes S^4H)=25$ and $\Pi_{25}$ is the associated projection map. In practice defining the projection maps $\Pi_k$ explicitly involves using various index symmetrisations, and contractions with $g$ and $\Omega$. For instance, $\Pi_{15}$ can be explicitly computed by  
\[
\begin{aligned}
(\Pi_{15}\alpha)_{i,jkl}
=\frac1{128}\Big[
&20g_{ij}\beta_{kl}
+5\bigl(\hat\Omega_{ij,kp}\beta_{pl}
        -\hat\Omega_{ij,lp}\beta_{pk}\bigr)
-3\hat\Omega_{kl,ip}\beta_{pj}\\
&-6\bigl(\hat\Omega_{jk,ip}\beta_{pl}
        -\hat\Omega_{jl,ip}\beta_{pk}\bigr)
-\bigl(\hat\Omega_{kl,jp}
       +2\hat\Omega_{jk,lp}
       -2\hat\Omega_{jl,kp}\bigr)\beta_{pi}
\Big],
\end{aligned}
\]
where $\beta_{ij}:=\alpha_{a,aij}\in \Lambda^2_{15}$ and $\hat{\Omega}:=\sum_{i=1}^3 \omega_i \otimes \omega_i$. Similar expressions hold for the other $\Pi_k$.

Recall from \eqref{decom symmetric tensor}, that any symmetric 2-tensor has 3 irreducible components. For $h \in S^2(\Lambda^1)$, we shall denote by $h_0$, $h_5$ and $h_{30}$ their irreducible components in view of \eqref{decom symmetric tensor}. In particular, we can apply this decomposition to the symmetric 2-tensors quadratic in torsion defined by $Y$ and $Z$ in \eqref{ioioioi}. The $7$ trivial representations in \eqref{quartic term} are generated by
\begin{equation}
|T|^4,\quad
|Y_{30}|^2,\quad
|Z_{30}|^2,\quad
| Y_{30}+Z_{30}|^2,\quad
|Y_5|^2,\quad
|(P_3\otimes P_3)B|^2,\quad
|(P_{10}\otimes P_{10})B|^2.
\end{equation}

Lastly we consider \eqref{qoqoqoq}. Recall that we already know that 2 generators are given by the divergence terms. Among the remaining 4 terms, 2 terms are eliminated by the Einstein condition i.e. $\Ric_0=0$, leaving only 2 unknown. A basis for these 2 terms are given by:

\begin{equation*}
\left\langle
\Pi_{70}\!\left((P_{15}\otimes P_{15})d^cT\right),
\Pi_{70}G
\right\rangle,\qquad 
\left\langle
\Pi_{25}\!\left((P_3\otimes P_{15})d^cT\right),
\Pi_{25}\!\left((P_3\otimes P_{15})B^t\right)
\right\rangle,
\end{equation*}
where $B^t_{abcd}:=B_{cdab}$, $G_{abcd}:=T_{i,ab}T_{i,cd}$, and $\Pi_k$ denote suitable $\Sp(2)\Sp(1)$ projection maps as before. 

Since we have now found the remaining 14 generators, we can apply the above strategy to extract the coefficient factors. This gives: 
\begin{align}
\mathcal H_{\mathrm{rem}}
={}&-4\operatorname{div}(\alpha_1)+4\operatorname{div}(\alpha_2)+\frac18|T|^4
+\frac{32}{15}|\Pi_{15}(\nabla^cT)|^2
 +\frac{32}{5}|\Pi_{25}(\nabla^cT)|^2\label{key identity}\\
 &+4|\Pi_{70}(\nabla^cT)|^2 +|Y_{30}|^2+4|Z_{30}|^2
 +|Y_5|^2
 +8|(P_3\otimes P_3)B|^2.\nonumber
\end{align}
Taking the trace in the last term $(P_3\otimes P_3)B \in S^2(\Lambda^2_3)\cong \R \oplus S^2_0(\Lambda^2_3)$ shows that 
\[
|(P_3\otimes P_3)B|^2=  \frac{1}{48}|T|^4+  |(P_3\otimes P_3)_0B|^2.
\]
Putting everything together, we can rewrite  \eqref{main integral identity 2} as
\[
\int_M \Big(\Big( \lambda+\frac{7}{24}|T|^2 \Big)|T|^2+H \Big)\vol=0,
\]
where $H$ is the sum of the remaining non-negative terms in \eqref{key identity}. This proves Theorem \ref{main thm}.

\bibliography{biblioG}
\end{document}

%% file: init.tex
\newcommand{\R}{\mathbb{R}}
\newcommand{\C}{\mathbb{C}}

\newcommand{\SO}{{\rm SO}}
\newcommand{\Sp}{{\rm Sp}}
\newcommand{\Spin}{{\rm Spin}}

\renewcommand{\epsilon}{\varepsilon}

\newcommand{\Ric}{{\rm Ric}}

\newcommand{\vol}{\mathrm{vol}}

\def\<{\mathopen{}\left<}
\def\>{\right>\mathclose{}}
\def\({\mathopen{}\left(}
\def\){\right)\mathclose{}}

\usepackage{multicol, color}

\definecolor{gold}{rgb}{0.85,.66,0}
\definecolor{cherry}{rgb}{0.9,.1,.2}
\definecolor{burgundy}{rgb}{0.8,.2,.2}
\definecolor{orangered}{rgb}{0.85,.3,0}
\definecolor{orange}{rgb}{0.85,.4,0}
\definecolor{olive}{rgb}{.45,.4,0}
\definecolor{lime}{rgb}{.6,.9,0}
\definecolor{green}{rgb}{.2,.7,0}
\definecolor{grey}{rgb}{.4,.4,.2}
\definecolor{brown}{rgb}{.4,.3,.1}

\newtheorem{theorem}{Theorem}
\newtheorem{proposition}{Proposition}

\theoremstyle{remark}
\newtheorem{remark}[equation]{Remark}

\theoremstyle{definition}

\def\w{\wedge}
\def\R{\mathbb{R}}
\def\C{\mathbb{C}}

\def\Lm{\Lambda}

\def\om{\omega}
\def\Om{\Omega}

\def\ip{\raise1pt\hbox{\large$\lrcorner\ \!$}} 